\documentclass[11pt,reqno]{amsart}
\usepackage[a4paper,left=1in,right=1in,top=1in,bottom=1in]{geometry}
\usepackage{amsmath,amssymb,amsthm,mathtools}
\usepackage{enumitem,cite}
\usepackage{tikz}
\usepackage{pgfplots}
\pgfplotsset{compat=1.18}
\usetikzlibrary{arrows.meta}
\usepackage[linkcolor=blue,citecolor=cyan,urlcolor=blue,colorlinks,hypertexnames=false]{hyperref}

\numberwithin{equation}{section}
\allowdisplaybreaks
\theoremstyle{plain}
\newtheorem{theorem}{Theorem}[section]
\newtheorem{proposition}[theorem]{Proposition}
\newtheorem{lemma}[theorem]{Lemma}

\theoremstyle{definition}

\newtheorem{remark}[theorem]{Remark}
\newcommand{\R}{\mathbb{R}}
\newcommand{\Z}{\mathbb{Z}}
\newcommand{\D}{\mathbb{D}}
\newcommand{\dd}{\mathrm d}
\newcommand{\dt}{\partial_t}
\newcommand{\Div}{\operatorname{div}}

\newcommand{\norm}[1]{\left\lVert#1\right\rVert}

\begin{document}
	\title[Ill-posedness of rotating free-boundary Euler]{%
		Rayleigh--Taylor instability and Hadamard ill-posedness
		for the rotating free-boundary Euler equations in three dimensions}
	
	\author{Chengchun Hao}
	\address{State Key Laboratory of Mathematical Sciences, Academy of Mathematics and Systems Science, Chinese Academy of Sciences, Beijing 100190, China \newline
		\and
		School of Mathematical Sciences, University of Chinese Academy of Sciences, Beijing 100049, China
	}
	\email{hcc@amss.ac.cn}
	\date{}

\begin{abstract}\sloppy
The free-boundary incompressible Euler equations linearized at a rigidly rotating equilibrium on the flat product of the unit disk and a circle are Hadamard ill-posed in Sobolev spaces of order at least two. The equilibrium violates the Taylor sign condition, under which local well-posedness is known. Eulerian linearization retaining free-surface advection reduces the spectral problem to a modified Bessel equation; Rouch\'e's theorem yields an unstable root whose real part grows as the square root of the axial frequency. Eigenmode Sobolev norms grow as the corresponding power of that frequency. Rescaling gives data tending to zero whose solutions are unbounded at any prescribed positive time. At zero axial frequency the analysis recovers Ebin's rotating-disk rate (Communications in Partial Differential Equations 12 (1987), 1175--1201). The product lift gives nonlinear ill-posedness on the invariant class of axially independent data.
\end{abstract}

\maketitle

\section{Introduction}\label{s:intro}

We study the free-boundary incompressible Euler equations on the flat product three-manifold
$\mathcal M:=\R^2\times S^1_L$, equipped with its product flat metric. The fluid occupies a moving
domain $\Omega_t\subset\mathcal M$ with free boundary $\Gamma_t:=\partial\Omega_t$ and unit outer
normal $N$. Write
$\mathcal G:=\{(t,x)\in[0,T]\times\mathcal M:x\in\Gamma_t\}$ for the space-time free boundary,
and set $D_t=\partial_t+u\cdot\nabla$. The divergence-free velocity $u$ and pressure $p$ satisfy
\begin{align}
	\partial_t u+u\cdot\nabla u+\nabla p&=0 &&\text{in }\Omega_t,\label{eq:euler}\\
	\operatorname{div}u&=0 &&\text{in }\Omega_t,\label{eq:div}\\
	D_t\text{ is tangent to }\mathcal G,\qquad p&=0 &&\text{on }\Gamma_t.
	\label{eq:bc}
\end{align}
All spatial derivatives are taken in the global parallel frame of the flat product metric. Under the Taylor sign condition $\nabla_Np\leqslant-c_0<0$ on the free surface for some constant $c_0>0$, local Sobolev well-posedness is known for Euclidean fluid domains. In the Euclidean state-space formulation of Ifrim--Pineau--Tataru--Taylor \cite{IPT}, the
threshold is $s>d/2+1$; in dimension three this is $s>5/2$, with the velocity and boundary charts
both measured at regularity $s$. The results below are proved directly on the flat product manifold
$\mathcal M$. Earlier a priori estimates and existence
results are given in \cite{CL,Lind05,CS,SZ08}; the irrotational case is due to Wu \cite{Wu97,Wu99}.

Hadamard-type ill-posedness for free-boundary Euler is known in specific settings. Ebin \cite{Ebin87}
proved discontinuous dependence at a rigidly rotating state in two dimensions, and Hao and Luo
\cite{HaoLuo2020} established ill-posedness for the ideal MHD free-boundary problem. Taylor
\cite{Taylor50} analyzed the instability of liquid surfaces accelerated normal to their planes;
Rayleigh \cite{Rayleigh83} studied the equilibrium of an incompressible heavy fluid of variable density.
Chandrasekhar's book \cite{Chandra} treats the classical hydrodynamic instability theory, including
rotation. Regarding the general necessity of the Taylor sign condition, Ifrim, Pineau, Tataru, and
Taylor describe Ebin's result as a specific ill-posedness scenario and state that they are unaware of
a proof in full generality \cite{IPT}. The results here concern the stated rotating state on a
flat product domain and do not give a general necessity theorem for the Taylor sign condition.

For the rotating unit disc, Ebin \cite{Ebin87} found linearized modes of azimuthal degree
$m\geqslant2$ with growth rate $\sqrt{m-1}$ and proved discontinuous dependence for the nonlinear
problem. We consider the same steady state on the flat product solid torus $\Omega=\D\times S^1_L$:
\[
u_0=x^\perp,\qquad p_0=\tfrac12(r^2-1),
\]
for which $\nabla_Np_0=1$ and the Taylor sign condition fails. The principal linear result here is an
explicit growing spectral root at nonzero axial frequency. For each fixed azimuthal mode
$m\geqslant0$ and all sufficiently large positive discrete axial frequencies
$k\in\frac{2\pi}{L}\mathbb Z$, the pressure equation reduces to a modified Bessel equation, and its
dispersion relation has a root satisfying
\[
\operatorname{Re}\lambda_k=\sqrt{k}\,(1+O(k^{-1})).
\]
The boundary-normalized eigenmodes have $H^{\mu_0}$ size comparable to $k^{\mu_0}$ for
$\mu_0\geqslant2$. Scaling them by rapidly decaying factors yields initial data tending to zero whose
linearized solutions are unbounded at each prescribed time $t_0>0$. The rate $\sqrt{|a|k}$, with
$|a|=1$, dominates the fixed angular frequency as $k\to\infty$.

At zero axial frequency, the growth is the two-dimensional rigid-rotation mechanism with rate
$\sqrt{m-1}$ for $m\geqslant2$. Separately, the nonlinear theorem below is the product lift of
Ebin's two-dimensional construction: it gives discontinuous dependence on the invariant class of
$x_3$-independent data on any positive common local existence interval. It does not assert a
nonlinear instability result for the nonzero axial-frequency modes; see
Theorem~\ref{thm:nonlin2d}.

Section~\ref{s:stat} does the steady state and its Taylor coefficient. Section~\ref{s:lin}
derives the Eulerian linearization and dispersion relation, including the two-dimensional anchor.
Section~\ref{s:3d} proves existence and asymptotics of the growing axial root, and
Section~\ref{s:thm} establishes unbounded growth of the linearized solution operator. The final
section proves the nonlinear result inherited from Ebin's construction.

Throughout, $r=(x_1^2+x_2^2)^{1/2}$, $\mathcal M=\R^2\times S^1_L$ with its flat product metric,
and $\Omega=\D\times S^1_L$, where $\D=\{r<1\}$ and $S^1_L=\R/L\Z$. The term \emph{solid torus}
refers to this product manifold and its topology; the model is not a curved Euclidean solid torus
embedded in $\R^3$. Since the base solution is stationary, we use the same symbol $\Omega$ for the
reference and current domains. The periodic coordinate $x_3$ and the global parallel frame identify
vectors with triples. We write $J=\begin{pmatrix}0&-1&0\\1&0&0\\0&0&0\end{pmatrix}$ and
$Jx=x^\perp=(-x_2,x_1,0)$; on scalar functions $J\cdot\nabla=\partial_\theta$. We use $y$ for
the label variable and $x$ for the physical variable. Sobolev norms are taken on the indicated
product domain. Finally, $f\asymp g$ means that their ratio is bounded above and below by positive
constants independent of the asymptotic parameter.

\section{The stationary solution}\label{s:stat}

\begin{proposition}[Stationary solution]
	\label{prop:stat}
	Let $L>0$ and $\Omega=\D\times S^1_L$. Then
	
	\begin{equation}\label{eq:steady}
		u_0(x)=Jx=x^\perp,\qquad p_0(x)=\tfrac12\big(r^2-1\big)
	\end{equation}
	is an exact solution of \eqref{eq:euler}--\eqref{eq:bc} on the fixed domain $\Omega$, with
	kinematic condition satisfied because $u_0\cdot N=0$ on $\{r=1\}$.
\end{proposition}

\begin{proof}
	$\Div(Jx)=0$. Since $J\cdot\nabla=\partial_\theta$ and $(Jx)\cdot N=(-x_2,x_1,0)\cdot(x_1,x_2,0)/r=0$
	on $\{r=1\}$, the boundary is stationary. Direct computation gives
	$\partial_1u_0=(0,1,0)$, $\partial_2u_0=(-1,0,0)$, hence
	$(u_0\cdot\nabla)u_0=-x_2(0,1,0)+x_1(-1,0,0)=(-x_1,-x_2,0)=-\nabla\big(\frac12r^2\big)$, so
	$\nabla p_0=-(u_0\cdot\nabla)u_0$. Finally $p_0=0$ on $\{r=1\}$.
\end{proof}

\begin{remark}[The Taylor coefficient]\label{rem:tt}
	Since $\nabla p_0=(x_1,x_2,0)$ and $N=(x_1,x_2,0)/r$, one has $\nabla_Np_0=r=1$ on $\{r=1\}$, so
	
	\begin{equation}\label{eq:a}
		a:=-\nabla_Np_0=-1<0
	\end{equation}
	on the whole free boundary. The Taylor sign condition is therefore violated, and it is violated
	non-degenerately. Equivalently, $p_0<0$ in the interior: the fluid is in tension, which is the state
	the sign condition excludes.
\end{remark}


\section{Linearization at the rotating steady state}\label{s:lin}

\subsection{The linearized system and the two boundary conditions}\label{ss:lin}

Let $\Omega=\D\times S^1_L$ and let \eqref{eq:steady} be the steady state of
Proposition~\ref{prop:stat}. Write $u=u_0+u^1$, $p=p_0+p^1$ and
$\Gamma_t=\{r=1+h(\theta,x_3,t)\}$. Linearizing \eqref{eq:euler}--\eqref{eq:bc} about
\eqref{eq:steady} gives

\begin{equation}\label{eq:lin}
	\dt u^1+(u_0\cdot\nabla)u^1+(u^1\cdot\nabla)u_0+\nabla p^1=0,\qquad
	\Div u^1=0\quad\text{in }\Omega ,
\end{equation}
and we look for modes

\begin{equation}\label{eq:mode0}
	\begin{gathered}
		u^1(r,\theta,x_3,t)=e^{\lambda t+im\theta+ikx_3}\,
		\big(u_r(r)e_r+u_\theta(r)e_\theta+u_z(r)e_3\big),\\
		p^1=e^{\lambda t+im\theta+ikx_3}\,p(r),\qquad m\in\Z,\ k\in\tfrac{2\pi}{L}\Z .
	\end{gathered}
\end{equation}
Since $u_0=x^\perp$ is a rigid rotation it is convenient to use the frame
$(e_r,e_\theta,e_3)$, with $e_3$ the unit vector along the $x_3$-axis, and the \emph{co-rotating}
frequency

\begin{equation}\label{eq:Lam}
	\Lambda:=\lambda+im\qquad\big(e^{\lambda t+im\theta}=e^{im\theta'}e^{\Lambda t},\ \ 
	\theta'=\theta-t\big).
\end{equation}

With $u^1$ as in \eqref{eq:mode0} and $(u^1\cdot\nabla)u_0=-u_\theta\,e_r+u_r\,e_\theta$, the
system \eqref{eq:lin} becomes, with $'=\frac{\dd}{\dd r}$,

\begin{equation}\label{eq:comp}
	\left\{
	\begin{aligned}
		\Lambda u_r-2u_\theta+p'&=0, &\qquad \Lambda u_\theta+2u_r+\tfrac{im}{r}p&=0,\\
		\Lambda u_z+ikp&=0,   &\qquad u_r'+\tfrac1r u_r+\tfrac{im}{r}u_\theta+iku_z&=0 .
	\end{aligned}
	\right.
\end{equation}
The boundary conditions linearize as follows. Since $\Gamma_t$ is a material surface and
$u_0\cdot N=0$ on $\{r=1\}$,

\begin{equation}\label{eq:kin}
	\dt h+(u_0\cdot\nabla)h=u^1\cdot N\quad\text{on }\{r=1\},
	\qquad\text{that is}\qquad \Lambda\,h=u_r(1);
\end{equation}
and since $p_0=0$ and $\partial_rp_0(1)=1$ on $\{r=1\}$,

\begin{equation}\label{eq:dyn}
	p^1+h\,\partial_rp_0=0\quad\text{on }\{r=1\},
	\qquad\text{that is}\qquad p(1)=-h .
\end{equation}
Eliminating $h$ gives the single spectral condition

\begin{equation}\label{eq:spec}
	\Lambda\,p(1)=-u_r(1)
\end{equation}
(for a steady state with Taylor coefficient $a=-\partial_rp_0(1)$ the right-hand side carries the
factor $a$; here $a=-1$, see \eqref{eq:a}).

\begin{remark}[The advection term in the kinematic condition]\label{rem:adv}
	The kinematic condition \eqref{eq:kin} reads $\Lambda h=u_r(1)$ rather than $\lambda h=u_r(1)$:
	the free surface is advected by the full velocity, so $u_0$ contributes $(u_0\cdot\nabla)h$ even
	when it is time independent. This is what yields the exact rate $\sqrt{m-1}$ at $k=0$
	(Theorem~\ref{thm:2d}).
\end{remark}

\subsection{Reduction to one equation for the pressure}\label{ss:red}

Assume $\Lambda\ne0$ and $\Lambda\neq\pm2i$, so that $A:=\Lambda^{2}+4$ is nonzero. No restriction
on $k$ is needed: the case $k=0$ reduces to the two-dimensional radial equation, and it is recovered below. The second and
third equations of \eqref{eq:comp} give
\[
u_\theta=\frac{-2u_r-\frac{im}{r}p}{\Lambda},\qquad u_z=-\frac{ikp}{\Lambda};
\]
substituting into the first gives

\begin{equation}\label{eq:usolve}
	u_r=-\frac{\Lambda p'+\frac{2im}{r}p}{A},\qquad A:=\Lambda^2+4 ,
\end{equation}
and substituting this, together with $u_\theta$ and $u_z$, into the fourth equation of \eqref{eq:comp},
and using $4-A=-\Lambda^2$, gives

\begin{equation}\label{eq:red}
	p''+\frac1r p'-\frac{m^2}{r^2}p=\frac{A\,k^2}{\Lambda^2}\,p.
\end{equation}
Thus the four functions reduce to a single second order equation, whose coefficient structure is
that of the two-dimensional Laplacian on the $e^{im\theta}$ mode, with the additional mass
$\mu^2:=k^2A/\Lambda^2$. Since $p$ is smooth at $r=0$,

\begin{equation}\label{eq:I}
	k=0:\quad p=r^{|m|};\qquad\qquad
	k\neq0:\quad p=I_m(\mu r),\qquad \mu=\frac{k\sqrt A}{\Lambda} ,
\end{equation}
where $I_m$ is the modified Bessel function and the branch of $\sqrt A$ is fixed by
$\mu\to k>0$ as $\lambda\to+\infty$ (see Section~\ref{ss:asym}). Substituting \eqref{eq:I} into
\eqref{eq:spec}, and using $p'(1)=\mu I_m'(\mu)$ together with
$\mu I_m'(\mu)=\mu I_{m-1}(\mu)-mI_m(\mu)$, we obtain the \emph{dispersion relation}

\begin{equation}\label{eq:DR}
	\Lambda A+m\Lambda-2im=\Lambda\mu\,\frac{I_{m-1}(\mu)}{I_m(\mu)},
	\qquad \mu=\frac{k\sqrt A}{\Lambda}.
\end{equation}
For the unstable modes studied below, $k>0$ and $\operatorname{Re}\Lambda>0$. Proposition~\ref{prop:mu}
then gives $\operatorname{Re}\mu>0$ and $I_m(\mu)\neq0$, so \eqref{eq:DR} is equivalent to
\eqref{eq:spec} after substituting \eqref{eq:usolve} and dividing by $I_m(\mu)$. The quotient form is
not used to characterize spectral values outside this half-plane. If $I_m(\mu)=0$ at a point with
$\Lambda\neq0,\pm2i$, then $\mu\neq0$. The Bessel equation is nonsingular at $\mu$, so uniqueness
for its initial-value problem implies $I_m'(\mu)\neq0$; otherwise $I_m$ would vanish identically.
Thus
\eqref{eq:usolve} gives
$u_r(1)=-\Lambda\mu I_m'(\mu)/A\neq0$, while \eqref{eq:spec} with $p(1)=0$ requires $u_r(1)=0$.
Thus such a point is not an eigenvalue. The case $k=0$ is treated in Section~\ref{ss:2d}, where the
ratio $I_{m-1}(\mu)/I_m(\mu)$ is not defined at $\mu=0$ for $m\geqslant1$.

\begin{remark}[Degenerate values of the spectral parameter]\label{rem:spurious}
	At $\Lambda=\pm2i$ one has $A=0$ and $\mu=0$, so the elimination \eqref{eq:usolve} is unavailable
	and \eqref{eq:comp} must be solved directly. For $\Lambda=2i$ the first two equations of \eqref{eq:comp} give $p'+mp/r=0$ for all $m\ge 0$. For
	$m=0$ one has $p\equiv c$, $u_\theta=iu_r$ and $u_z=-kc/2$. The divergence equation gives
	$(ru_r)'=ik^2cr/2$, hence $u_r=ik^2cr/4+C/r$; regularity at $r=0$ forces $C=0$. The spectral
	condition \eqref{eq:spec} then becomes $c(1+k^2/8)=0$, so $c=0$ and the mode is trivial. For
	$m\geqslant1$, the regularity of $p$ at $r=0$ forces $p=cr^{-m}$ to vanish. Then $u_\theta=iu_r$,
	$u_z=0$, and the divergence equation gives $u_r=C r^{m-1}$. The spectral condition requires $u_r(1)=0$, so $C=0$. At $\Lambda=-2i$ the first
	two equations similarly imply, for $m=0$, $u_\theta=-iu_r$ and $p\equiv c$; the divergence
	equation gives $u_r=-ik^2cr/4+C/r$, with $C=0$ by regularity. The spectral condition again forces
	$c=0$. For $m\geqslant1$, consistency gives $p=cr^m$, and the fourth equation admits the unique
	regular solution
	\[
	u_r(r)=-\frac{icm}{4}r^{m-1}-\frac{ick^2}{4(m+1)}r^{m+1},
	\]
	whose spectral condition $-2ic=\frac{ic}{4}\bigl(m+\frac{k^2}{m+1}\bigr)$ forces $c=0$.
	At $\Lambda=0$ one has $p\equiv0$ for $k\neq0$, and the rest of \eqref{eq:comp} gives the trivial
	solution. Hence for $k\neq0$ none of $0,\pm2i$ is an eigenvalue for any $m$, and in the right
	half-plane the eigenvalues are exactly the roots of \eqref{eq:DR}. For $k=0$ and $m=0$, however,
	$\Lambda=0$ is a genuine neutral mode ($u_r=u_\theta=u_z=0$, $p\equiv1$, $h=-1$), which the
	reduction to \eqref{eq:DR} misses because it divides by $\Lambda$: it corresponds to the
	infinitesimal change of radius in the family of stationary rigidly rotating discs.
\end{remark}

\begin{proposition}[The Bessel parameter]\label{prop:mu}
	Let $k>0$, let $\mathbb H:=\{\Lambda:\operatorname{Re}\Lambda>0\}$, and let
	\begin{equation}\label{eq:mu-def}
		\mu(\Lambda):=k\sqrt{1+\frac4{\Lambda^2}},\qquad \Lambda\in\mathbb H ,
	\end{equation}
	the branch being the one fixed by $\mu\to k$ as $\Lambda\to+\infty$, i.e.\ the branch of
	\eqref{eq:I}. Then
	\begin{enumerate}
		\item[(i)] $\mu$ is holomorphic on $\mathbb H$ and $\operatorname{Re}\mu(\Lambda)>0$ there;
		\item[(ii)] $\operatorname{Im}\mu(\Lambda)$ has the opposite sign to $\operatorname{Im}\Lambda$ (and
		$\operatorname{Im}\mu=0$ if and only if $\operatorname{Im}\Lambda=0$);
		\item[(iii)] $I_m(\mu(\Lambda))\neq0$ for every $m\in\mathbb N_0$ and every $\Lambda\in\mathbb H$;
		\item[(iv)] the branch points $\Lambda=\pm2i$ and all the points at which $I_m(\mu)=0$ lie on the
		imaginary axis.
	\end{enumerate}
\end{proposition}

\begin{proof}
	\emph{Holomorphy and the sign of the real part.} If $1+4/\Lambda^2$ belonged to $(-\infty,0]$ then
	$\Lambda^2$ would be real and non-positive, so that $\Lambda\in i\mathbb R$, which is excluded by
	$\operatorname{Re}\Lambda>0$. Hence $\Lambda\mapsto1+4/\Lambda^2$ sends $\mathbb H$ into
	$\mathbb C\setminus(-\infty,0]$, on which the principal square root is holomorphic and has strictly
	positive real part; the composition is therefore holomorphic on $\mathbb H$, is positive on the
	positive real axis, and satisfies $\operatorname{Re}\mu>0$. This is the branch of \eqref{eq:mu-def},
	since $\mu\to k$ as $\Lambda\to+\infty$.
	
	\emph{The sign of the imaginary part.} Write $\Lambda=a+ib$ with $a>0$. Then
	\[
	\Im\frac4{\Lambda^2}=\Im\frac{4\overline{\Lambda}^{\,2}}{|\Lambda|^{4}}
	=-\frac{8ab}{|\Lambda|^{4}},
	\]
	so that $\Im(1+4/\Lambda^2)$ has the sign of $-b$. Writing $1+4/\Lambda^2=\rho e^{i\theta}$ with
	$\theta\in(-\pi,\pi)$, the principal square root is $\sqrt\rho\,e^{i\theta/2}$, so
	$\Im\mu=k\sqrt\rho\,\sin(\theta/2)$ has the sign of $\theta$, hence the sign of $-b$, as claimed.
	
	\emph{The zeros of $I_m$.} Since $I_m(z)=i^{-m}J_m(iz)$ and all zeros of integer-order $J_m$ are real \cite[10.21(i)]{DLMF}, $I_m(z)=0$
	forces $z$ to be purely imaginary, which is excluded by $\operatorname{Re}\mu>0$.
	
	\emph{Position of the singular points.} For $m\geqslant1$, $I_m(\mu)=0$ gives $\mu=0$ or
	$\mu^{2}=-\beta^{2}$ with $\beta>0$ a zero of $J_m$; in the first case $1+4/\Lambda^2=0$, i.e.\
	$\Lambda=\pm2i$, and in the second
	\[
	\Lambda^{2}=\frac{4}{-\beta^{2}/k^{2}-1}=-\frac{4k^{2}}{\beta^{2}+k^{2}}<0 ,
	\]
	so that $\Lambda\in i\mathbb R$. For $m=0$ the same conclusion holds, $I_0(0)=1\neq0$.
\end{proof}

The equations \eqref{eq:euler}--\eqref{eq:bc} and the steady state are real; the steady state is
invariant under $x_3\mapsto-x_3$ but anti-invariant under $\theta\mapsto-\theta$ (it changes sign
there). Hence conjugating a mode $(m,k)$ produces a mode $(-m,-k)$ with the same growth rate
$\operatorname{Re}\lambda$ (its $\lambda$ becomes $\bar\lambda$), and reflecting $x_3$ produces
$(m,-k)$. It suffices to treat $m\geqslant0$ and $k>0$, which is what we do below.

We have normalized the rotation rate to $1$; for a base state $u_0=\omega x^\perp$ with general
$\omega>0$ the same computation goes through with
\[
\Lambda=\lambda+im\omega,\qquad A=\Lambda^2+4\omega^2,\qquad
\mu=\frac{k\sqrt A}{\Lambda}=k\sqrt{1+\frac{4\omega^2}{\Lambda^2}} ,
\]
and with \eqref{eq:DR} replaced by
$\Lambda A+\omega^2m\Lambda-2i\omega^3m=\omega^2\Lambda\mu\,I_{m-1}(\mu)/I_m(\mu)$: the factor
$\omega^2$ in front of $m\Lambda$ and of $\Lambda\mu$ comes from $\partial_rp_0(1)=\omega^2$, which
is also what multiplies the spectral condition \eqref{eq:spec}. Then
$\operatorname{Re}\lambda\approx\omega\sqrt k=\sqrt{|a|k}$, the Rayleigh--Taylor rate with effective
gravity $|a|=\omega^2$ (the leading balance is
$\Lambda^3\approx\omega^2\Lambda\mu\approx\omega^2\Lambda k$). Only the
\emph{ratio} $I_{m-1}/I_m$ occurs in \eqref{eq:DR}, the individual functions being of size $e^\mu$.

\subsection{The two-dimensional anchor}\label{ss:2d}

For $k=0$, \eqref{eq:red} is the homogeneous Euler equation, whose regular solution is
$p=r^{|m|}$; hence $p'(1)/p(1)=|m|$ and \eqref{eq:spec} becomes, for $m\geqslant0$,

\begin{equation}\label{eq:anchor}
	\Lambda^3+(4-m)\Lambda-2im=0 .
\end{equation}
The root $\Lambda=-2i$ is spurious, as explained in Remark~\ref{rem:spurious}. For
$m\geqslant2$, the genuine growing root is $\Lambda=i+\sqrt{m-1}$. Since $\lambda=\Lambda-im$,
it gives

\begin{equation}\label{eq:anchor2}
	\lambda=\sqrt{m-1}-i(m-1)\qquad(m\geqslant2).
\end{equation}

\begin{theorem}\label{thm:2d}
	For every $m\geqslant2$ the linearized problem \eqref{eq:lin} with \eqref{eq:kin}--\eqref{eq:dyn} has a solution of the
	form \eqref{eq:mode0} with $k=0$ and $\lambda$ as in \eqref{eq:anchor2}; in particular
	$\operatorname{Re}\lambda=\sqrt{m-1}$.
\end{theorem}

\begin{proof}
	For $k=0$ the reduced equation \eqref{eq:red} becomes $p''+\frac1r p'-\frac{m^2}{r^2}p=0$,
	whose solution smooth at the origin is $p(r)=r^m$ (constant normalized to $1$). With
	$\Lambda=i+\sqrt{m-1}$ (the growing root of \eqref{eq:anchor}; for $m\geqslant2$ we have
	$A:=\Lambda^2+4\neq0$), define
	\[
	u_r(r)=-\frac{\Lambda p'(r)+\frac{2im}{r}p(r)}{A},\qquad
	u_\theta(r)=\frac{-2u_r(r)-\frac{im}{r}p(r)}{\Lambda},\qquad
	u_z(r)=0,
	\]
	which by the reduction in \S\ref{ss:red} satisfy all four equations of \eqref{eq:comp}. For the
	boundary conditions, \eqref{eq:dyn} is satisfied by taking $h:=-p(1)=-1$; and \eqref{eq:kin}
	holds because at $r=1$
	\[
	u_r(1)=-\frac{m}{A}(\Lambda+2i)=-\Lambda=\Lambda h,
	\]
	the middle equality being the identity $A\Lambda=m(\Lambda+2i)$ that follows from the anchor
	polynomial \eqref{eq:anchor}, since
	$(\Lambda^2+4)\Lambda=\Lambda^3+4\Lambda=(m-4)\Lambda+2im+4\Lambda=m(\Lambda+2i)$. Thus
	\eqref{eq:mode0} with this data and $\lambda=\Lambda-im=\sqrt{m-1}-i(m-1)$ is the asserted
	solution.
\end{proof}

\section{The genuinely three-dimensional modes: \texorpdfstring{$k\neq0$}{k nonzero}}\label{s:3d}

\subsection{Asymptotics: growth rate \texorpdfstring{$\sqrt k$}{sqrt k}}\label{ss:asym}

Fix $m$ and let $k\gg1$; the following heuristic is made rigorous in
Section~\ref{ss:root}. Suppose that $\lambda$ is of order $\sqrt k$, so that
$\Lambda\approx\lambda$, $A=\Lambda^2+4\approx\Lambda^2$ and

\begin{equation}\label{eq:muasy}
	\mu=\frac{k\sqrt A}{\Lambda}\approx k .
\end{equation}
For fixed $m$, Lemma~\ref{lem:bessel} below gives
$I_{m-1}(\mu)/I_m(\mu)=1+\frac{2m-1}{2\mu}+O(\mu^{-2})$. Under the assumed scale
$\Lambda^2\asymp k$, the definition of $\mu$ gives $\mu=k+O(1)$. Dividing \eqref{eq:DR} by
$\Lambda$ and using this Bessel expansion yields
\[
\Lambda^2+4+m-\frac{2im}{\Lambda}
=\mu+m-\frac12+O(\mu^{-1}),
\qquad\text{hence}\qquad \Lambda^2=k+O(1),
\]
that is

\begin{equation}\label{eq:main}
	\lambda=\sqrt k-im+O(k^{-1/2}),\qquad
	\operatorname{Re}\lambda=\sqrt k\,\big(1+O(k^{-1})\big),\quad
	\operatorname{Im}\lambda=-m+O(k^{-1/2}).
\end{equation}
For $m=0$ one correction term can be computed: here \eqref{eq:spec} reads
$(\lambda^2+4)=\mu I_1(\mu)/I_0(\mu)$ with $\mu=k\sqrt{1+4/\lambda^2}$, and inserting
$I_1/I_0=1-\frac1{2\mu}-\frac1{8\mu^2}+O(\mu^{-3})$, that is
$\mu I_1(\mu)/I_0(\mu)=\mu-\frac12-\frac1{8\mu}+O(\mu^{-2})$, and
$\mu=k+\frac{2k}{\lambda^2}+O(k^{-1})=k+2+O(k^{-1})$, gives
$\lambda^2+4=k+\frac32+O(k^{-1})$, i.e.

\begin{equation}\label{eq:corr}
	\lambda^2=k-\tfrac52+O(k^{-1})\qquad(m=0).
\end{equation}

\begin{remark}[The modal fields are not horizontal]\label{rem:horiz}
	By the third equation of \eqref{eq:comp}, $u_z=-ikp/\Lambda\neq0$ for $k\neq0$, so every nontrivial
	separated mode has a nonzero third component and $\partial_3p^1=ikp\neq0$. In particular, there is
	no nontrivial $x_3$-dependent horizontal growing mode in the separated class \eqref{eq:mode0}: if
	$u_z\equiv0$, the third equation gives $p\equiv0$, and the first two give $(\Lambda^2+4)u_r=0$, so
	$u_r=u_\theta=0$ for a growing mode with $\Lambda^2+4\neq0$.
\end{remark}

\subsection{Existence and asymptotics of the root}\label{ss:root}

The asymptotics of Section~\ref{ss:asym} were formal. Here the existence of the root, and its
asymptotics, are proved for every fixed $m$ and all large $k$. Multiplying \eqref{eq:DR} by
$\Lambda\ne0$ and
using $A=\Lambda^2+4$ and $I_{-1}=I_1$, the dispersion relation is equivalent to $F(\Lambda)=0$,
where

\begin{equation}\label{eq:F}
	F(\Lambda):=\Lambda^4+(4+m)\Lambda^2-2im\Lambda-\Lambda^2\,\mu\,\frac{I_{m-1}(\mu)}{I_m(\mu)},
	\qquad \mu=k\sqrt{1+\frac4{\Lambda^2}} ,
\end{equation}
the branch of the root being fixed by $\mu\to k$ as $\Lambda\to+\infty$ (this is the form
\eqref{eq:I} of $\mu$, since $\mu=k\sqrt A/\Lambda=k\sqrt{1+4/\Lambda^2}$). Only the \emph{ratio}
$I_{m-1}/I_m$ is transcendental.

\begin{lemma}[The Bessel ratio]\label{lem:bessel}
	Let $m\geqslant0$ be an integer. One has $I_m(\mu)\ne0$ whenever $\operatorname{Re}\mu>0$.
	As $|\mu|\to\infty$ in any fixed closed subsector
	$|\arg\mu|\leqslant\frac\pi2-\delta$, with $\delta>0$, one has, uniformly,
	\begin{equation}\label{eq:Iratio}
		\mu\,\frac{I_{m-1}(\mu)}{I_m(\mu)}
		=\mu+\Big(m-\tfrac12\Big)+\frac{\alpha_m}{\mu}
		+O\big(|\mu|^{-2}\big),\qquad \alpha_m:=\frac{4m^2-1}{8}.
	\end{equation}
	The restriction to a closed subsector of the right half-plane is essential; no expansion across the
	imaginary axis is asserted.
	Finally, with $r_m(\mu):=I_{m-1}(\mu)/I_m(\mu)$ there is the \emph{exact} identity
	\begin{equation}\label{eq:fratio}
		\frac{\dd}{\dd\mu}\Big(\mu\,\frac{I_{m-1}(\mu)}{I_m(\mu)}\Big)
		=\mu\Big(1-r_m(\mu)^{2}\Big)+2m\,r_m(\mu),
	\end{equation}
	and in particular, for $m=0$,
	\begin{equation}\label{eq:f0}
		f_0'(\mu)=\mu\Big(1-\Big(\frac{I_1(\mu)}{I_0(\mu)}\Big)^{2}\Big)>0
		\qquad(\mu>0),
	\end{equation}
	so that $f_0(\mu):=\mu I_1(\mu)/I_0(\mu)$ is strictly increasing on $(0,\infty)$ and tends to
	$+\infty$ there. The inequalities $0<I_1(\mu)<I_0(\mu)$ for $\mu>0$ follow from the series and from
	$I_0(\mu)-I_1(\mu)=\frac1\pi\int_0^\pi e^{\mu\cos\theta}(1-\cos\theta)\,\dd\theta>0$.
\end{lemma}

\begin{proof}
	$I_m(\mu)=i^{-m}J_m(i\mu)$ and all zeros of integer-order $J_m$ are real \cite[10.21(i)]{DLMF}, so $I_m(\mu)=0$ forces $i\mu\in\R$,
	i.e.\ $\mu$ purely imaginary; this is excluded when $\operatorname{Re}\mu>0$.
	
	For $|\arg\mu|\leqslant\frac\pi2-\delta$, the fixed-order expansion
	\[
	I_\nu(\mu)=\frac{e^\mu}{\sqrt{2\pi\mu}}
	\left(1-\frac{4\nu^2-1}{8\mu}
	+\frac{(4\nu^2-1)(4\nu^2-9)}{128\mu^2}+O(\mu^{-3})\right)
	\]
	is uniform \cite[10.40.1]{DLMF}. Dividing the expansions for orders $m-1$ and $m$ gives
	\[
	\frac{I_{m-1}}{I_m}
	=1+\frac{\alpha_m-\alpha_{m-1}}{\mu}
	+\frac{\alpha_m^2-\alpha_{m-1}\alpha_m+\beta_{m-1}-\beta_m}{\mu^2}
	+O(\mu^{-3}),
	\]
	where $\alpha_\nu=\frac{4\nu^2-1}{8}$ and
	$\beta_\nu=\frac{(4\nu^2-1)(4\nu^2-9)}{128}$. The first coefficient is
	$\alpha_m-\alpha_{m-1}=m-\frac12$, and the second is $\alpha_m$, since
	\[
	\alpha_m^2-\alpha_{m-1}\alpha_m+\beta_{m-1}-\beta_m
	=\alpha_m\Big(m-\tfrac12\Big)+\alpha_m\Big(\tfrac32-m\Big)=\alpha_m,
	\]
	using $\beta_{m-1}-\beta_m=\alpha_m(\frac32-m)$. Multiplication by $\mu$ proves
	\eqref{eq:Iratio}.
	
	Note the sign of the leading correction: it is $+\big(m-\frac12\big)$, so that, for instance,
	$\mu I_1(\mu)/I_0(\mu)=\mu-\frac12-\frac{1}{8\mu}+O(\mu^{-2})$ for $\operatorname{Re}\mu>0$; the sign
	is the one carried by the expansion above.
	
	For \eqref{eq:fratio} one uses $I_0'=I_1$, $I_1'=I_0-I_1/\mu$ and the recurrence
	$I_{m-2}-I_m=\frac{2(m-1)}{\mu}I_{m-1}$ (with $I_{-1}=I_1$, $I_{-2}=I_2$), which give
	\[
	\frac{\dd}{\dd\mu}\Big(\frac{I_{m-1}}{I_m}\Big)
	=1-\Big(\frac{I_{m-1}}{I_m}\Big)^{2}+\frac{2m-1}{\mu}\frac{I_{m-1}}{I_m} ,
	\]
	and \eqref{eq:fratio} follows on multiplying by $\mu$ and adding $r_m$. For $m=0$ this reads
	$f_0'=\mu(1-r_0^{2})$, which is positive because $0<r_0<1$ on $(0,\infty)$; the two inequalities
	follow from the series $I_1=\sum_{k\geqslant0}\frac{(\mu/2)^{2k+1}}{k!(k+1)!}>0$ and from the
	integral representation quoted above.
\end{proof}

\begin{lemma}[Rescaled convergence]\label{lem:resc}
	Put $\Lambda=\sqrt k\,\zeta$. Then, uniformly for $\zeta$ in compact subsets of
	$\{\operatorname{Re}\zeta>0\}$,
	
	\begin{equation}\label{eq:muk}
		\mu=k+\frac2{\zeta^2}-\frac2{\zeta^4k}+O\big(k^{-2}\big),
	\end{equation}
	and, with
	
	\begin{equation}\label{eq:Psi}
		\Psi_k(\zeta):=k^{-2}F\big(\sqrt k\,\zeta\big),\qquad
		\Psi_\infty(\zeta):=\zeta^4-\zeta^2=\zeta^2(\zeta-1)(\zeta+1) ,
	\end{equation}
	one has, with $\alpha_m=(4m^2-1)/8$ as in \eqref{eq:Iratio},
	
	\begin{equation}\label{eq:Psiasy}
		\Psi_k(\zeta)=\Psi_\infty(\zeta)+\frac1k\Big[\tfrac92\,\zeta^2-2\Big]
		-2im\zeta\,k^{-3/2}+\frac{2/\zeta^2-\alpha_m\zeta^2}{k^2}+O\big(k^{-3}\big) .
	\end{equation}
\end{lemma}

\begin{proof}
	$\mu=k\sqrt{1+4/(k\zeta^2)}=k\big(1+\frac2{k\zeta^2}-\frac2{k^2\zeta^4}+O(k^{-3})\big)
	=k+\frac2{\zeta^2}-\frac2{\zeta^4k}+O(k^{-2})$, which is \eqref{eq:muk}. Dividing \eqref{eq:F} by
	$k^2$ with $\Lambda=\sqrt k\zeta$ gives $\zeta^4$, $(4+m)\zeta^2/k$, $-2im\zeta k^{-3/2}$ and
	$-\zeta^2\cdot\frac{\mu}k\frac{I_{m-1}}{I_m}$; by \eqref{eq:Iratio} and \eqref{eq:muk},
	\[
	\begin{aligned}
		\frac{\mu}k\frac{I_{m-1}}{I_m}
		&=\frac{k+\frac2{\zeta^2}-\frac2{\zeta^4k}+\big(m-\frac12\big)+\frac{\alpha_m}{k}+O(k^{-2})}{k}\\
		&=1+\frac{\frac2{\zeta^2}+\big(m-\frac12\big)}{k}
		+\frac{\alpha_m-\frac2{\zeta^4}}{k^2}+O(k^{-3}) ,
	\end{aligned}
	\]
	so the last term is
	$-\zeta^2-\frac{2+\big(m-\frac12\big)\zeta^2}{k}+\frac{2/\zeta^2-\alpha_m\zeta^2}{k^2}+O(k^{-3})$.
	Adding the four contributions gives \eqref{eq:Psiasy}, since
	$(4+m)-\big(m-\frac12\big)=\frac92$. Of the four, only $-2im\zeta k^{-3/2}$ is of order
	$k^{-3/2}$; it is absent when $m=0$. At $\zeta=1$ the expansion reduces to
	$\Psi_k(1)=\frac{5}{2k}-\frac{2im}{k^{3/2}}+\frac{2-\alpha_m}{k^2}+O(k^{-3})$, which matches
	the expansion of \eqref{eq:F} directly.
\end{proof}

\begin{theorem}[Existence and asymptotics]\label{thm:root}
	Let $m\geqslant0$ be fixed. There are $k_0=k_0(m)$, $C=C(m)$ and $\rho=\rho(m)\in(0,\frac18)$ such that
	for every $k\geqslant k_0$ the dispersion relation \eqref{eq:DR} has a root $\Lambda_k$, unique after the
	rescaling $\Lambda=\sqrt k\,\zeta$ in the disc $|\zeta-1|<\rho$, with
	$\Lambda_k=\sqrt k+O(k^{-1/2})$. More precisely, with
	$c_m:=\frac{8m^2+21}{32}$,
	
	\begin{equation}\label{eq:rootsharp}
		\Lambda_k=\sqrt k-\frac{5}{4\sqrt k}+\frac{im}{k}+\frac{c_m}{k^{3/2}}+\frac{im}{2k^{2}}
		+O\big(k^{-5/2}\big),
	\end{equation}
	so that the real and imaginary parts are
	$\operatorname{Re}\Lambda_k=\sqrt k-\frac{5}{4\sqrt k}+\frac{c_m}{k^{3/2}}+O(k^{-5/2})$ and
	$\operatorname{Im}\Lambda_k=\frac{m}{k}+\frac{m}{2k^{2}}+O(k^{-5/2})=\frac{m}{k}+O(k^{-2})$. In particular the coefficient of $k^{-1/2}$
	does not depend on $m$. Consequently
	
	\begin{equation}\label{eq:rootasym}
		\lambda_k=\Lambda_k-im,\qquad
		\operatorname{Re}\lambda_k=\sqrt k\,\big(1+O(k^{-1})\big),\qquad
		\operatorname{Im}\lambda_k=-m+O\big(k^{-1}\big) ,
	\end{equation}
	and, for $m=0$, $\lambda_k^2=k-\frac52+O(k^{-1})$.
\end{theorem}

\begin{proof}
	Let $\rho\in(0,\frac18)$. On $|\zeta-1|=\rho$ we have $|\zeta|\geqslant1-\rho\geqslant\frac78$ and
	$|\zeta+1|\geqslant2-\rho\geqslant\frac{15}8$, hence
	$|\Psi_\infty|=|\zeta|^2|\zeta-1||\zeta+1|\geqslant(\frac78)^2\rho\frac{15}8=:c_\rho>0$, while
	$\sup_{|\zeta-1|\leqslant\rho}|\Psi_k-\Psi_\infty|\leqslant C_\rho/k\to0$ by Lemma~\ref{lem:resc},
	the estimate on the closed disc following from the one on its boundary by the maximum modulus principle. Rouch\'e's
	theorem gives, for $k\geqslant k_0$, that $\Psi_k$ and $\Psi_\infty$ have equally many zeros in
	$D(1,\rho)$, i.e.\ exactly one: $\Psi_\infty=\zeta^2(\zeta-1)(\zeta+1)$ has there the simple zero
	$1$ only (the zeros $-1$ and the double zero $0$ lie outside, as $\rho<\frac18$). Let $\zeta_k$ be
	it. For $m=0$, the rescaled equation has real coefficients on the conjugation-invariant disk
	$D(1,\rho)$, so $\Psi_k(\overline\zeta)=\overline{\Psi_k(\zeta)}$ there. Uniqueness of the zero
	in this disk implies $\overline{\zeta_k}=\zeta_k$, and hence $\Lambda_k$ is real. Since
	$\Psi_\infty=\zeta^2(\zeta-1)(\zeta+1)$ and, on $D(1,\rho)$ with
	$\rho<\frac18$, both $|\zeta|^2\geqslant\frac{49}{64}$ and $|\zeta+1|\geqslant\frac{15}8$, one has
	$|\Psi_\infty(\zeta)|\geqslant\frac{735}{512}\,|\zeta-1|$ there. Hence
	$|\zeta_k-1|\leqslant\frac{512}{735}\,|\Psi_\infty(\zeta_k)|
	=|\Psi_k(\zeta_k)-\Psi_\infty(\zeta_k)|\cdot\frac{512}{735}\leqslant C/k$, whence
	$\zeta_k=1+O(k^{-1})$ and $\Lambda_k=\sqrt k\,\zeta_k=\sqrt k+O(k^{-1/2})$. Finally
	$\operatorname{Re}\lambda_k=\sqrt k\,\operatorname{Re}\zeta_k=\sqrt k(1+O(k^{-1}))$ and
	$\operatorname{Im}\lambda_k=\sqrt k\,\Im\zeta_k-m=-m+O(k^{-1/2})$.
	
	To justify the sharper expansion \eqref{eq:rootsharp}, set $\varepsilon=k^{-1/2}$ and write
	$\eta=\zeta-1$. The uniform expansion \eqref{eq:Psiasy} on a fixed neighborhood of $1$ becomes
	\[
	\begin{aligned}
		\Psi_k(\zeta)&=\Psi_\infty(\zeta)+\varepsilon^2B(\zeta)
		-2im\zeta\varepsilon^3+\varepsilon^4C(\zeta)+O(\varepsilon^6),\\
		B(\zeta)&=\tfrac92\zeta^2-2,\qquad
		C(\zeta)=2/\zeta^2-\alpha_m\zeta^2,
	\end{aligned}
	\]
	with a remainder uniform in $\zeta$. Consider the truncated approximation
	\[
	\eta_* = a\varepsilon^2+ib\varepsilon^3+c\varepsilon^4+id\varepsilon^5.
	\]
	Substitution into the displayed equation and cancellation of the coefficients of
	$\varepsilon^2,\varepsilon^3,\varepsilon^4,\varepsilon^5$, respectively, gives
	\[
	\begin{gathered}
		2a+\tfrac52=0,\qquad 2b-2m=0,\\
		2c+5a^2+9a+(2-\alpha_m)=0,\qquad
		2d+10ab+9b-2ma=0.
	\end{gathered}
	\]
	Thus $a=-\frac54$, $b=m$, $c=\frac{8m^2+21}{32}$, and $d=\frac m2$. The uniform remainder
	in \eqref{eq:Psiasy} then gives $|\Psi_k(1+\eta_*)|\leqslant C\varepsilon^6$.
	
	It remains to compare this approximate zero with the Rouch\'e zero $\zeta_k$. The functions
	$\Psi_k$ are holomorphic on a neighborhood of the closed disk under consideration and converge
	locally uniformly to $\Psi_\infty$; Cauchy's formula therefore gives locally uniform convergence of
	their first derivatives. Since $\Psi_\infty'(1)=2$, after reducing the disk if necessary, the
	factorization
	\[
	\Psi_k(\zeta)=(\zeta-\zeta_k)G_k(\zeta),\qquad
	G_k(\zeta)=\int_0^1\Psi_k'\big(\zeta_k+t(\zeta-\zeta_k)\big)\,\dd t,
	\]
	has $|G_k(\zeta)|\geqslant c>0$ for $\zeta$ in that disk and all sufficiently large $k$. Hence
	$|\zeta_k-(1+\eta_*)|\leqslant C\varepsilon^6$. Retaining the terms through order $\varepsilon^5$ in
	$\zeta_k=1+\eta_*$ leaves a remainder of order $\varepsilon^6$; since
	$\Lambda_k=\sqrt{k}\,\zeta_k=\varepsilon^{-1}\zeta_k$, it follows that
	\[
	\Lambda_k=\sqrt k-\frac{5}{4\sqrt k}+\frac{im}{k}
	+\frac{8m^2+21}{32k^{3/2}}+\frac{im}{2k^2}+O(k^{-5/2}),
	\]
	which is \eqref{eq:rootsharp}. For $m=0$, squaring gives
	$\Lambda_k^2=k-\frac52+\frac{23}{8k}+O(k^{-2})$, and thus
	$\lambda_k^2=k-\frac52+O(k^{-1})$, as claimed.
\end{proof}

\begin{remark}[The surface components of the eigenmode]\label{rem:surfcomp}
	The mode of Theorem~\ref{thm:root} satisfies, at the free surface, $u_r(1)=-\Lambda_k$ by
	\eqref{eq:spec} and $u_z(1)=-ik/\Lambda_k$ by the third equation of \eqref{eq:comp}, with the
	normalization $p(1)=1$ of Lemma~\ref{lem:size}. Hence
	\[
	\frac{|u_z(1)|}{|u_r(1)|}=\frac{k}{|\Lambda_k|^{2}}=1+O\big(k^{-1}\big),
	\]
	the last equality following from \eqref{eq:rootsharp}: the axial and the radial surface velocities
	are of the same size.
\end{remark}

\begin{remark}[Two mechanisms at $k=0$ and at $k\neq0$]\label{rem:two}
	At $k=0$ the Rayleigh--Taylor scaling vanishes, and the growth in Theorem~\ref{thm:2d} comes
	instead from the rigid-rotation mechanism of \cite{Ebin87}, with rate $\sqrt{m-1}$ for $m\geqslant2$.
	The two asymptotic families are the two-dimensional modes $(k=0$, $m\geqslant2)$ and the axial modes
	$(m$ fixed, $k\to\infty)$, which are distinct limits of the dispersion relation \eqref{eq:DR}.
\end{remark}

\begin{remark}[Localization of the eigenvectors]\label{rem:bl}
	Since $\mu\approx k$, the profile $p=I_m(\mu r)$ varies on the scale $1/k$: the eigenvectors are
	concentrated in a boundary layer of width $1/k$ at the free surface, and
	$p'(1)/p(1)\approx\mu\approx k$. A thin boundary layer makes the normal derivative of the pressure
	at the surface large, and \eqref{eq:spec} converts it into $\Lambda^2\approx k$. Consequently the
	normalized eigenprofiles develop derivatives that grow with $k$. The layer variable is $\mu_k(1-r)$, and since $\mu_k=k+2+O(k^{-1})$ the width $1/\mu_k$ is $1/k$ to leading order.
\end{remark}


\section{Ill-posedness of the linearized problem}\label{s:thm}

The reference free boundary is the flat torus $\Gamma_0:=\{r=1\}$, so that the
mode data of Lemma~\ref{lem:size} have boundary displacement in
$H^{\mu_0+1/2}(\Gamma_0)$ and the nonlinear state space measures it in
$H^s(\Gamma_0)$.

\begin{lemma}[Size of the eigenmode]\label{lem:size}
	Fix $m\geqslant0$ and $L>0$, let $\mu_0\geqslant2$ be real (the non-integral case is treated at
	the end of the proof), and let $k$ be a discrete axial frequency
	$k\in\frac{2\pi}{L}\Z$ with $k\geqslant k_0(m)$, $k_0$ as in Theorem~\ref{thm:root}. For the root
	$\lambda_k$ of \eqref{eq:DR} put $\mu_k=k\sqrt A/\Lambda_k$ and
	\[
	p^{(k)}(r):=\frac{I_m(\mu_kr)}{I_m(\mu_k)},
	\]
	so that $p^{(k)}(1)=1$; let $u^{(k)}$ be the associated field of
	\eqref{eq:comp}--\eqref{eq:usolve}. There are positive constants $c_{\mu_0,m,L}$ and
	$C_{\mu_0,m,L}$ such that
	
	\begin{equation}\label{eq:size}
		c_{\mu_0,m,L}k^{\mu_0}\leqslant \norm{u^{(k)}}_{H^{\mu_0}(\Omega)}\leqslant C_{\mu_0,m,L}k^{\mu_0}.
	\end{equation}
	Thus $\norm{u^{(k)}}_{H^{\mu_0}(\Omega)}\asymp k^{\mu_0}$, with implied constants depending on
	$\mu_0,m,L$.
	In the state norm of Theorem~\ref{thm:lin}, the boundary displacement belongs to
	$H^{\mu_0+1/2}(\Gamma_0)$. Since $\dim\Gamma_0=2$ and $\mu_0+\frac12\geqslant\frac52>2$,
	Sobolev embedding gives $h\in C^1(\Gamma_0)$, so the free surface is a $C^1$ graph and the
	kinematic/dynamic boundary conditions \eqref{eq:kin}--\eqref{eq:dyn} are well defined.
\end{lemma}

\begin{proof}
	All the fields of \eqref{eq:comp}--\eqref{eq:usolve} are of the form
	$e^{im\theta+ikx_3}F(r)$, with $F$ an algebraic combination of $I_m$, its derivatives, powers of
	$r$ and the constants $\Lambda,A$; here $\mu=\mu_k\asymp k$, while $\operatorname{Re}\Lambda=\operatorname{Re}\lambda\asymp\sqrt k$ and
	$|\operatorname{Im}\Lambda|=O(1)$ by \eqref{eq:rootasym}, so that $\Lambda=O(k^{1/2})$, $A\asymp k$
	and $A^{-1}=O(k^{-1})$.
	
	Write $M=|\mu|$. In the sector $|\arg\mu|\leqslant\frac\pi2-\delta$, one has
	$\operatorname{Re}\mu\geqslant c_\delta M$. For each fixed integer order $\ell\geqslant0$ and
	integer $j\geqslant0$, the integral representation
	\[
	I_\ell(z)=\frac1\pi\int_0^\pi e^{z\cos\theta}\cos(\ell\theta)\,\dd\theta
	\]
	and its $j$ derivatives give
	\[
	|I_\ell^{(j)}(\mu r)|\leqslant C_{\ell,j}
	\frac{e^{r\operatorname{Re}\mu}}{(1+Mr)^{1/2}},\qquad 0\leqslant r\leqslant1.
	\]
	The large-argument expansion in this sector also gives
	$|I_m(\mu)|\geqslant c_m e^{\operatorname{Re}\mu}M^{-1/2}$ for sufficiently large $M$. Therefore, for
	$\phi(r)=I_m(\mu r)/I_m(\mu)$,
	\[
	|\partial_r^j\phi(r)|\leqslant C_{m,j}M^j
	\left(\frac{M}{1+Mr}\right)^{1/2}e^{-c_\delta M(1-r)}.
	\]
	Splitting the radial integral at $r=1/2$ shows
	$\int_0^1\frac{Mr}{1+Mr}e^{-2c_\delta M(1-r)}\,\dd r\leqslant C M^{-1}$; hence
	
	\begin{equation}\label{eq:bessellayer}
		\int_0^1\left|\left(\frac1\mu\frac{\dd}{\dd r}\right)^{j}\phi(r)\right|^2r\,\dd r
		\leqslant C_{m,j}M^{-1},\qquad
		\int_0^1|\partial_r^j\phi(r)|^2r\,\dd r\leqslant C_{m,j}M^{2j-1}.
	\end{equation}
	The constants are uniform for $m$ in any fixed finite set and for $\mu$ in the stated sector. This
	argument controls the region near the axis as well as the boundary layer; it does not use a boundary
	asymptotic outside its range of validity.
	
	For Cartesian derivatives, apply
	$\partial_1+i\partial_2$ and $\partial_1-i\partial_2$ to the full regular mode
	$I_m(\mu r)e^{im\theta}$. The Bessel recurrences give
	\[
	\mu I_{m+1}(\mu r)e^{i(m+1)\theta}
	\quad\text{and}\quad
	\mu I_{m-1}(\mu r)e^{i(m-1)\theta},
	\]
	respectively (with $I_{-q}=I_q$). Thus every fixed-order Cartesian derivative
	is a finite sum of normalized profiles of fixed neighboring orders, multiplied by the corresponding
	power of $\mu$. The same integral estimate applies to each such order. The factors $r^{-1}$ in
	\eqref{eq:usolve} and in $u_\theta$ are eliminated using
	$I_{m-1}(z)-I_{m+1}(z)=2mI_m(z)/z$, so the regularity at $r=0$ is included in these bounds.
	
	Using the recurrences above in \eqref{eq:usolve} and the formulae for $u_\theta,u_z$, each
	component is a finite sum of normalized neighboring Bessel profiles with coefficients bounded by
	$C_m k^{1/2}$ (since $\mu_k\asymp k$, $\Lambda_k\asymp k^{1/2}$, and $A\asymp k$). By
	\eqref{eq:bessellayer}, each profile has $L^2(r\,\dd r)$ norm at most $C_m k^{-1/2}$, giving
	$\norm{u^{(k)}}_{L^2(\Omega)}\leqslant C_{m,L}$. For the lower bound, on the layer
	$1-|\mu_k|^{-1}\leqslant r\leqslant1$ the large-argument expansion gives
	$|I_m(\mu_k r)/I_m(\mu_k)|\geqslant c_m$. Since
	$u_z^{(k)}=-ikp^{(k)}/\Lambda_k$ and $k/|\Lambda_k|\asymp k^{1/2}$, integration over this layer,
	whose radial measure is comparable to $k^{-1}$, yields
	$\norm{u^{(k)}}_{L^2(\Omega)}\geqslant c_{m,L}$. Thus
	$\norm{u^{(k)}}_{L^2(\Omega)}\asymp1$. An axial derivative contributes a factor $k$; a
	Cartesian derivative of order $j$ contributes at most $C_{m,j}k^j$ by the neighboring-order
	Bessel bounds. Consequently
	$\norm{D^\beta u^{(k)}}_{L^2(\Omega)}\leqslant C_{m,L,\beta}k^{|\beta|}$
	for every fixed multi-index $\beta$, and in particular for $|\beta|\leqslant\mu_0$, whence
	$\norm{u^{(k)}}_{H^{\mu_0}}\leqslant C_{\mu_0,m,L}k^{\mu_0}$. Conversely, the $\mu_0$ axial derivatives carry the factor
	$e^{ikx_3}$ and contribute $\norm{D_3^{\mu_0}u^{(k)}}_{L^2(\Omega)}\asymp k^{\mu_0}
	\norm{u^{(k)}}_{L^2(\Omega)}\asymp k^{\mu_0}$, so that
	$\norm{u^{(k)}}_{H^{\mu_0}}\geqslant c_{\mu_0,m,L}k^{\mu_0}$, proving \eqref{eq:size} for integer
	$\mu_0$.
	
	\emph{Nonintegral Sobolev orders.} We use the restriction Sobolev spaces on the smooth product
	domain $\Omega=\D\times S^1_L$. They satisfy the standard complex interpolation identity
	\[
	(H^{\ell_0}(\Omega),H^{\ell_1}(\Omega))_{[\theta]}=H^{\mu_0}(\Omega),
	\qquad \mu_0=(1-\theta)\ell_0+\theta\ell_1.
	\]
	Applying this to the integer-order upper bounds above gives
	\[
	\|u^{(k)}\|_{H^{\mu_0}(\Omega)}\leqslant C_{\mu_0,m,L}k^{\mu_0}\|u^{(k)}\|_{L^2(\Omega)}.
	\]
	For the lower bound, the axial Fourier projection is contractive on these product Sobolev spaces,
	and on the single mode $e^{ikx_3}$ the axial Bessel-potential multiplier has symbol
	$(1+k^2)^{\mu_0/2}$. Consequently
	\[
	\|u^{(k)}\|_{H^{\mu_0}(\Omega)}
	\geqslant (1+k^2)^{\mu_0/2}\|u^{(k)}\|_{L^2(\Omega)}.
	\]
	Thus $\|u^{(k)}\|_{H^{\mu_0}}\asymp k^{\mu_0}$ for real $\mu_0\geqslant0$, which proves
	\eqref{eq:size}.
\end{proof}

\begin{theorem}[Unbounded growth of the linearized solution operator]\label{thm:lin}
	Fix an azimuthal mode $m\geqslant0$. Let $L>0$, $\Omega=\D\times S^1_L$, $t_0>0$ and $\mu_0\geqslant2$. Then there is a sequence of initial data for the linearized problem \eqref{eq:lin}--\eqref{eq:dyn}, indexed by the discrete axial frequencies $k_n=\frac{2\pi n}{L}\to\infty$, and constants depending on $L,\mu_0,m$, such that
	\[
	\norm{u^{(n)}(0)}_{H^{\mu_0}(\Omega)}+\norm{h^{(n)}(0)}_{H^{\mu_0+\frac12}(\{r=1\})}
	\leqslant C_{\mu_0,L,m}e^{-\frac12\operatorname{Re}\lambda_{k_n} t_0}
	\]
	and
	\[
	\norm{u^{(n)}(t_0)}_{H^{\mu_0}(\Omega)}\ \geqslant\ c_{\mu_0,L,m}k_n^{-1}e^{+\frac12\operatorname{Re}\lambda_{k_n} t_0}
	\ \xrightarrow[n\to\infty]{}\ \infty ,
	\]
	where $\operatorname{Re}\lambda_{k_n}=\sqrt{k_n}\,(1+O(k_n^{-1}))$ for the eigenvalue of
	\eqref{eq:DR} in Theorem~\ref{thm:root} (applied with $k=k_n$).
	Taking real and imaginary parts of the complex modes (or, equivalently, combining the conjugate axial frequencies $k_n$ and $-k_n$) gives real-valued data with the same estimates. In particular the data tend to $0$ in $H^{\mu_0}$ while the corresponding solutions are unbounded in $H^{\mu_0}$: no constant $C$ bounds $\|u(t_0)\|_{H^{\mu_0}}$ by
	$C\|u(0)\|_{H^{\mu_0}}$ on the modes of the family, and the linearized problem is not Hadamard well-posed.
\end{theorem}

\begin{proof}
	Fix $m\geqslant0$. Theorem~\ref{thm:root} provides the eigenvalue $\lambda_{k_n}$, with $\operatorname{Re}\lambda_{k_n}=\sqrt{k_n}\,(1+O(k_n^{-1}))$, for every fixed $m$ and
	all large $k_n$. Let $u^{(n)},p^{(n)}$ be
	the eigenmode of Lemma~\ref{lem:size} with
	$k=k_n$. Set $M:=\mu_0+1$ and
	\[
	u^{(n)}(0):=A_n\,u^{(n)},\qquad h^{(n)}(0):=A_n\,u^{(n)}_r(1)/\Lambda_{k_n} ,
	\qquad A_n:=k_n^{-M}e^{-\frac12\operatorname{Re}\lambda_{k_n}t_0} .
	\]
	By \eqref{eq:kin} the corresponding solution of the linearized problem is
	$A_ne^{\lambda_{k_n}t}$ times the mode, so, by the asymptotic $\|u^{(k_n)}\|_{H^{\mu_0}}\asymp k_n^{\mu_0}$
	in Lemma~\ref{lem:size},
	\[
	\norm{u^{(n)}(t_0)}\ \geqslant\ A_nc_{\mu_0,L,m}k_n^{\mu_0}e^{\operatorname{Re}\lambda_{k_n}t_0}
	\ =\ c_{\mu_0,L,m}k_n^{-1}e^{\frac12\operatorname{Re}\lambda_{k_n}t_0}\ \longrightarrow\ \infty ,
	\]
	whereas $\norm{u^{(n)}(0)}\leqslant C_{\mu_0,L,m}A_nk_n^{\mu_0}=C_{\mu_0,L,m}k_n^{-1}e^{-\frac12\operatorname{Re}\lambda_{k_n}t_0}\to0$.
	For the surface displacement, combining \eqref{eq:spec} with \eqref{eq:kin} gives at once
	$h^{(n)}=u^{(n)}_r(1)/\Lambda_{k_n}=-p^{(n)}(1)$,
	so with the normalization $p^{(n)}(1)=1$ of Lemma~\ref{lem:size}
	\[
	\begin{aligned}
		h^{(n)}(0)&=-A_ne^{im\theta+ik_nx_3},\qquad\text{hence}\\
		\norm{h^{(n)}(0)}_{H^{\mu_0+\frac12}(\Gamma_0)}
		&\leqslant C_{\mu_0,L,m}k_n^{\mu_0+\frac12}A_n
		\leqslant C_{\mu_0,L,m}k_n^{M}A_n ,
	\end{aligned}
	\]
	which is the bound required in the statement. (Here $|u^{(n)}_r(1)|=|\Lambda_{k_n}|\asymp\sqrt{k_n}$
	is \emph{not} $O(1)$; the two factors of $\sqrt{k_n}$ cancel in
	$h^{(n)}(0)=u^{(n)}_r(1)/\Lambda_{k_n}$, so $h^{(n)}$ has the same angular and axial frequencies as
	the normal velocity, as required by \eqref{eq:kin}.)
\end{proof}


\section{Nonlinear ill-posedness inherited from the two-dimensional problem}\label{s:nonlin2d}

\begin{theorem}[Nonlinear ill-posedness for $x_3$-independent data]\label{thm:nonlin2d}
	Let $s\geqslant4$ be an integer and $\Omega_h=\{r<1+h(\theta,x_3)\}$. Fix a smooth cutoff
	$\chi_0(r)$ supported in a collar of $r=1$ and equal to $1$ at $r=1$, and define the flattening
	\[
	\Psi_h(r,\theta,x_3):=(r+\chi_0(r)h(\theta,x_3),\theta,x_3),
	\qquad \Psi_0=\mathrm{Id},
	\]
	for $h$ in a sufficiently small $H^s(\Gamma_0)$ neighborhood of zero. This family and its inverses
	are uniformly bounded as $H^s$ composition operators. Define
	\[
	\|(u,h)-(u_0,0)\|_{\mathcal X^s}
	:=\|u\circ\Psi_h-u_0\|_{H^s(\Omega)}+\|h\|_{H^s(\Gamma_0)}.
	\]
	Let $\mathcal X^s_{\mathrm{ax}}$ be the relative state space of data independent of $x_3$. No $T>0$
	is admissible: the solution map restricted to this invariant class is not Hadamard well-posed on a
	relative $\mathcal X^s$-neighborhood of $(u_0,0)$, with continuous dependence into
	$C([0,T];\mathcal X^s_{\mathrm{ax}})$.
\end{theorem}

\begin{proof}
	Suppose the problem is locally Hadamard well-posed on a relative $\mathcal X^s$-neighborhood of $(u_0,0)$ in $\mathcal X^s_{\mathrm{ax}}$: solutions exist,
	are unique, and depend continuously on the data as a map into $C([0,T];\mathcal X^s_{\mathrm{ax}})$.
	Fix $t_0\in(0,T]$. The assumed solution map supplies existence and uniqueness on the common interval
	$[0,T]$, which are the hypotheses used in Ebin's Section 5 argument \cite{Ebin87}. Applying that
	argument gives smooth disk data converging to the rotating disk state, and corresponding unique disk
	flows $\eta_n^{(2)}$ with boundary displacement $h_n$. The base rotation is the flow $R_t=\exp(tJ)$ of
	the steady velocity $u_0$: $\partial_tR_t=u_0\circ R_t=J\,R_t$ and $R_0=\mathrm{Id}$. In terms of this
	flow,
	\[
	\|(\eta_n^{(2)}(t_0)-R_{t_0},\partial_t\eta_n^{(2)}(t_0)-\partial_tR_{t_0})\|_{H^1\times H^1}
	\longrightarrow\infty.
	\]
	Ebin's initial configurations are the identity and the velocity perturbation is
	$e^{-n^{1/4}}z^n$, where $z=x_1+ix_2$ is the complex planar coordinate on the disk $\D$; hence, for
	each fixed $s$, its $H^s(\D)$ norm tends to zero, since derivatives of $z^n$ grow only
	polynomially in $n$. With $h_n(0)=0$ and $\Psi_0=\mathrm{Id}$, these data therefore converge in
	$\mathcal X^s$. The assumed local well-posedness places them in its neighborhood and, by uniqueness,
	identifies Ebin's constructed trajectories with those of the solution map. Extend each solution to
	the flat product solid torus by
	\[
	\begin{aligned}
		\eta_n(y,x_3,t)&:=(\eta_n^{(2)}(y,t),x_3),\\
		u_n(x_1,x_2,x_3,t)&:=(u_n^{(2)}(x_1,x_2,t),0),\\
		\Omega_{n,t}&:=\eta_n(t,\Omega),
	\end{aligned}
	\]
	Because the fields and boundary are independent of $x_3$, this is a solution of the Euler equations
	on $\mathcal M=\R^2\times S^1_L$. The initial-data norms on the disk and the product domain differ
	only by fixed factors depending on $L$.
	
	Write
	\[
	w_n(t):=u_n(t)\circ\Psi_{h_n(t)}-u_0,
	\qquad v_n(t):=u_n(t)-u_0\quad\text{on }\Omega_{n,t}.
	\]
	Continuity into $C([0,T];\mathcal X^s_{\mathrm{ax}})$ gives
	$\sup_{t\in[0,T]}(\|w_n(t)\|_{H^s(\Omega)}+\|h_n(t)\|_{H^s(\Gamma_0)})\to0$.
	By the specified radial collar map, $\|\Psi_{h_n(t)}-\mathrm{Id}\|_{H^s(\Omega)}\leqslant C\|h_n(t)\|_{H^s(\Gamma_0)}$.
	Since $u_0(x)=Jx$ in the planar components,
	\[
	v_n(t)\circ\Psi_{h_n(t)}
	=w_n(t)+u_0-u_0\circ\Psi_{h_n(t)},
	\]
	so the uniform composition bounds for $\Psi_h$ and $\Psi_h^{-1}$ yield
	\[
	\sup_{t\in[0,T]}\|v_n(t)\|_{H^s(\Omega_{n,t})}
	\leqslant C\sup_{t\in[0,T]}\big(\|w_n(t)\|_{H^s(\Omega)}+\|h_n(t)\|_{H^s(\Gamma_0)}\big)
	\longrightarrow0.
	\]
	Since $h_n\to0$ in $C([0,T];H^s(\Gamma_0))$, discard finitely many indices so that the global
	extensions $\Phi_{n,t}$ of the collar maps $\Psi_{h_n(t)}$, equal to the identity outside a fixed
	collar, are diffeomorphisms of $\mathcal M$ with uniformly bounded $H^s$ composition norms for
	$\Phi_{n,t}$ and $\Phi_{n,t}^{-1}$. They satisfy $\Phi_{n,t}(\Omega)=\Omega_{n,t}$. Fix a bounded
	linear extension operator $E_0:H^s(\Omega)\to H^s(\mathcal M)$, supported in a fixed bounded
	neighborhood of $\Omega$. For $f\in H^s(\Omega_{n,t})$, define
	\[
	E_{n,t}f:=\bigl(E_0(f\circ\Phi_{n,t})\bigr)\circ\Phi_{n,t}^{-1}.
	\]
	Because $E_0$ is a right inverse to restriction to $\Omega$, $E_{n,t}f=f$ on $\Omega_{n,t}$.
	Uniform composition bounds give
	\[
	\|E_{n,t}f\|_{H^s(\mathcal M)}\leqslant C\|f\|_{H^s(\Omega_{n,t})},
	\]
	with $C$ independent of $n,t$; the extensions have support in a common bounded set and depend
	continuously on $t$ when applied to the continuous paths $v_n(t)$.
	
	Choose a smooth cutoff $\chi$ supported in a fixed bounded set and equal to one on a neighborhood
	of every $\Omega_{n,t}$. Define the extended velocity
	\[
	\widehat u_n(t):=\chi u_0+E_{n,t}v_n(t).
	\]
	It agrees with $u_n$ on $\Omega_{n,t}$ and satisfies
	$\sup_{n,t\in[0,T]}\|\widehat u_n(t)\|_{H^s(\mathcal M)}\leqslant C_T$ by the preceding extension
	bound and $\sup_t\|v_n(t)\|_{H^s(\Omega_{n,t})}\to0$. Let $\widehat\eta_n$ be its flow from the
	identity on $\mathcal M$. Since $\widehat u_n=u_n$ along the moving fluid domain, uniqueness for the particle ODE gives
	$\widehat\eta_n|_\Omega=\eta_n$. The fields $\widehat u_n$ have common compact support, so restrict
	them to a fixed smooth compact subdomain containing that support and equal to zero near its boundary.
	For $s\geqslant4>3/2+1$, the standard Sobolev flow estimate on this fixed domain
	\cite{EM70} gives, whenever $\|W\|_{C([0,T];H^s)}\leqslant M_0$,
	\[
	\sup_{t\in[0,T]}\bigl(\|\eta_W(t)-\mathrm{Id}\|_{H^s}
	+\|\eta_W(t)^{-1}-\mathrm{Id}\|_{H^s}\bigr)\leqslant C(s,T,M_0),
	\]
	where $\eta_W$ is the flow of $W$. Indeed, $\partial_tD\eta_W=(DW)\circ\eta_W\,D\eta_W$ and
	$H^s\hookrightarrow W^{1,\infty}$ give uniform bounds for $D\eta_W$ and its inverse. For
	$2\leqslant q\leqslant s$, differentiating the flow equation $q$ times yields the leading term
	$(DW)\circ\eta_W\,D^q\eta_W$; every other term is a product of derivatives of $W$ of order at most
	$q$ composed with $\eta_W$ and derivatives of $\eta_W$ of order at most $q-1$. The change-of-variables
	formula, the already established first-derivative bounds, and the Sobolev product estimates bound
	these terms in $L^2$ by $C(s,M_0)$ times the lower-order flow norms. Induction on $q$ and Gronwall give
	the stated $H^s$ bound. The inverse flow satisfies
	$\partial_t\eta_W^{-1}+(D\eta_W^{-1})W=0$ and the same estimates apply. Applying this estimate with
	$W=\widehat u_n$ and the uniform bound above gives
	\[
	\sup_{n,t\in[0,T]}\bigl(\|\widehat\eta_n(t)-\mathrm{Id}\|_{H^s(\mathcal M)}
	+\|\widehat\eta_n(t)^{-1}-\mathrm{Id}\|_{H^s(\mathcal M)}\bigr)\leqslant C(s,T,M_0).
	\]
	On $\Omega$, the extension property and the standard Sobolev composition estimate give
	\[
	\|v_n(t)\circ\eta_n(t)\|_{H^s(\Omega)}
	=\|(E_{n,t}v_n(t))\circ\widehat\eta_n(t)\|_{H^s(\Omega)}
	\leqslant C_T\|v_n(t)\|_{H^s(\Omega_{n,t})},
	\]
	where $C_T$ is independent of $n,t$.
	
	Let $\eta_0(t)=R_t$ be the steady flow and let $z_n=\eta_n-\eta_0$, interpreted in the global
	parallel frame; the third components agree identically. The initial domains agree, so $z_n(0)=0$.
	Since $u_0(x)=Jx$ in the planar components, the flow equations give
	\[
	\partial_tz_n=Jz_n+v_n\circ\eta_n.
	\]
	The composition estimate and Gronwall's inequality imply
	\[
	\sup_{0\leqslant t\leqslant t_0}\|z_n(t)\|_{H^s(\Omega)}
	\leqslant C_{t_0}\sup_{0\leqslant t\leqslant t_0}
	\|v_n(t)\|_{H^s(\Omega_{n,t})}\longrightarrow0.
	\]
	The flow identity then gives
	$\sup_{0\leqslant t\leqslant t_0}\|\partial_tz_n(t)\|_{H^s(\Omega)}\to0$ as well. This
	contradicts Ebin's $H^1$ divergence at $t=t_0$, since $H^s(\D)$ embeds continuously in $H^1(\D)$.
	Hence no such continuous local solution map exists.
\end{proof}

\section*{Acknowledgments}
Hao's research was partially supported by the National Natural Science Foundation of China (NSFC) under Grant No. 12171460, the CAS Project for Young Scientists in Basic Research under Grant No. YSBR-031 and the National Key R\&D Program of China under Grant No. 2021YFA1000800.

\section*{AI assistance statement} The author used AI tools to assist with exploring, numerical verification and manuscript editing; all mathematical validation, final proof decisions, and final wording remain the sole responsibility of the author.

\end{document}